\documentclass[12pt]{article}
\usepackage{amssymb}

\usepackage{amsmath,amsfonts}
\usepackage{amscd}

\newtheorem{lemma}{Lemma}[section]
\newtheorem{coro}[lemma]{Corollary}

\newtheorem{thm}[lemma]{Theorem}
\newtheorem{defn}[lemma]{Definition}
\makeatletter\@addtoreset{equation}{section}
\renewcommand\theequation{\thesection.\@arabic\c@equation}

\begin{document}
\begin{center}
{\LARGE   Generalized Obata theorem and its applications on
foliations}

 \renewcommand{\thefootnote}{}
 \footnote {This paper was supported by
KRF-2008-313-C00076 from Korea Research Foundations and
R01-2008-000-20370-0 from KOSEF.}
\renewcommand{\thefootnote}{\arabic{footnote}}
\setcounter{footnote}{0}

\vspace{1 cm} {\large Seoung Dal Jung, Keum Ran Lee and Ken
Richardson}
\end{center}
\vspace{0.5cm}

{\bf Abstract.} We prove the generalized Obata theorem on
foliations. Let $M$ be a complete Riemannian manifold with a
foliation $\mathcal F$ of codimension $q\geq 2$ and a bundle-like
metric $g_M$. Then $(M,\mathcal F)$ is transversally isometric to
$(S^q(1/c),G)$, where $S^q(1/c)$ is the $q$-sphere of radius $1/c$
in $(q+1)$-dimensional Euclidean space and $G$ is a discrete
subgroup of the orthogonal group $O(q)$, if and only if there
exists a non-constant basic function $f$ such that $\nabla_X
df=-c^2 fX^b$ for all basic normal vector fields $X$, where $c$ is
a positive constant and $\nabla$ is the connection on the normal
bundle. By the generalized Obata theorem, we classify such
manifolds which admit transversal non-isometric conformal fields.

\section{Introduction}
In 1962, M. Obata [\ref{Obata}] proved that a complete Riemannian
manifold $(M,g_M)$ is isometric with a sphere  of radius $\frac1c$
in (n+1)-dimensional Euclidean space if and only if $M$ admits a
non-constant function $f$ such that
\begin{align*}
\nabla^{M}_X df = -c^2 fX^b
\end{align*}
for any vector $X$, where $\nabla^M$ is a Levi-Civita connection
on $M$, $c$ is positive constant and $X^b$ is the $g_M$-dual form
of $X$.

In Section 2, we recall basic facts on foliations as well as some
lemmas we need for the main results. In Section 3, we generalize
the Obata theorem to Riemannian foliations. Namely, let $M$ be a
complete Riemannian manifold with a foliation $\mathcal F$ of
codimension $q\geq 2$ and a bundle-like metric $g_M$. Then the
foliation $(M,\mathcal F)$ is transversally isometric to
$(S^q(1/c),G)$, where $S^q(1/c)$ is the $q$-sphere of radius $1/c$
in the $(q+1)$-Euclidean space and $G$ is a discrete subgroup of
$O(q)$, if and only if there exists a non-constant basic function
$f$ such that $\nabla_X df=-c^2 fX^b$ for all basic normal vector
fields $X$, where $\nabla$ is the connection on the normal bundle.

 A Riemannian foliation is a foliation
$\mathcal F$ on a smooth $n$-manifold $M$ such that the quotient
bundle $N\mathcal F\cong Q=TM/T\mathcal F$ is endowed with a
metric $g_Q$ satisfying $\theta(X)g_Q=0$ for any vector $X\in
T\mathcal F$, where $T\mathcal F$ is the tangent bundle of
$\mathcal F$ and $\theta(X)$ is the transverse Lie derivative
([\ref{Tond}]). Note that we can choose a Riemannian metric $g_M$
on $M$ such that $g_{N\mathcal{F}}:=g_M|_{N\mathcal{F}} =g_Q$;
such a metric is called {\it bundle-like}. In the last section, as
applications of the generalized Obata theorem,  we study the
Riemannian foliations admitting transversal non-isometric
conformal fields.
\section{Preliminaries}
Let $(M,g_M,\mathcal F)$ be a $(p+q)$-dimensional Riemannian
manifold with a foliation $\mathcal F$ of codimension $q$ and a
bundle-like metric $g_M$ with respect to $\mathcal F$.
 Then we have an exact sequence of vector bundles
\begin{align}\label{eq1-1}
 0 \longrightarrow T\mathcal F \longrightarrow
TM {\overset\pi\longrightarrow} N\mathcal F \longrightarrow 0,
\end{align}
where $T\mathcal F$ is the tangent bundle  and $N\mathcal F\cong
Q=TM/T\mathcal F$ is the normal bundle of $\mathcal F$. We denote
by $\nabla$ the connection on the normal bundle
$N\mathcal{F}\subset TM$. That is, $\nabla_X Y=\pi(\nabla^{M}_X
Y)$ for any $X,Y\in N\mathcal F$, where $\nabla^{M}$ is the
Levi-Civita connection on $M$; this connection $\nabla$ is
guaranteed to be metric and torsion-free with respect to
$g_Q=g_{N\mathcal{F}}$ (\cite{{Tond},{Tond1}}). Let $R^\nabla,
K^\nabla,\rho^\nabla$ and $\sigma^\nabla$ be the transversal
curvature tensor, transversal sectional curvature, transversal
Ricci operator and transversal scalar curvature with respect to
$\nabla$, respectively. A differential form $\omega\in
\Omega^r(M)$ is {\it basic} if
\begin{equation}
i(X)\omega=0,\ \theta(X)\omega=0, \quad \forall X\in T\mathcal F.
\end{equation}
Let $\Omega_B^r(\mathcal F)$ be the set of all basic r-forms on
$M$. Then $L^2(\Omega^*(M))$ is decomposed as [\ref{Lop},
\ref{Park}]
\begin{align}\label{eq1-5}
L^2(\Omega(M))=L^2(\Omega_B(\mathcal F)) \oplus
L^2(\Omega_B(\mathcal F))^\perp.
\end{align}
Let $P:L^2(\Omega^*(M))\to L^2(\Omega_B^*(\mathcal F))$ be the
orthogonal projection onto basic forms [\ref{Park}], which
preserves smoothness in the case of Riemannian foliations. For any
$r$-form $\phi$, we put the basic part of $\phi$ as $\phi_B
:=P\phi$.
 Let $\delta_B$ be the formal adjoint operator of
$d_B=d|_{\Omega_B^*(\mathcal F)}$.  The basic Laplacian $\Delta_B$
is given by
 $ \Delta_B = d_B\delta_B
+ \delta_B d_B$. For any basic function $f$, it is well-known
(\cite{jj} or \cite{Park}) that $\int_M \Delta_B f =0$.
\begin{lemma}  Let $(M,g_M,\mathcal F)$ be a closed, connected Riemannian manifold with a foliation $\mathcal
F$ and a bundle-like metric $g_M$. If $\Delta_B f =
\kappa_B^\sharp(f)$ for any basic function $f$, then $f$ is
constant.
\end{lemma}
{\bf Proof.} Since $\Delta_B$ on basic functions is the
restriction of the elliptic operator
$\Delta+\left(\kappa_B^\sharp-\kappa^\sharp\right)$ on all
functions (see \cite[Prop. 4.1]{Park}), we have that
$\Delta_B-\kappa_B^\sharp$ on basic functions is the restriction
of $\Delta-\kappa^\sharp$ on all functions. A solution $u:M\to
\mathbb{R}$ to $\left(\Delta-\kappa^\sharp\right)u=0$ satisfies
the maximum and minimum principles locally, so that if such a $u$
has a local maximum or minimum, then $u$ is constant. The result
follows.  $\Box$

Let $V(\mathcal F)$ be the space of all vector fields $Y$ on $M$
satisfying $[Y,Z]\in T\mathcal F$ for all $Z\in T\mathcal F$. An
element of $V(\mathcal F)$ is called an {\it infinitesimal
automorphism} of $\mathcal F$. Let
\begin{align}\label{eq1-2}
\overline V(\mathcal F)=\{\overline Y:=\pi(Y)\ |\ Y\in V(\mathcal
F)\}.
\end{align}
Then we have $\Omega_B^r(\mathcal F)\subset \Gamma(\Lambda^r Q^*)$
and $\overline V(\mathcal F)\cong \Omega_B^1(\mathcal F)$.
 If $Y\in
V(\mathcal F)$ satisfies $\theta(Y)g_Q=0$, then $\overline Y$ is
called a {\it transversal Killing field} of $\mathcal F$.
 If $Y\in V(\mathcal F)$ satisfies $\theta(Y)g_Q=2f_Y g_Q$ for
a basic function $f_Y$ depending on $Y$, then $\overline Y$ is
called a {\it transversal conformal field} of $\mathcal F$; in
this case, we have
\begin{equation}
f_Y = \frac 1q \operatorname{div_\nabla} \overline Y.
\end{equation}
A $\overline Y$ is called a {\it transversal non-isometric
conformal field} of $\mathcal F$ if $\overline Y$ is transversal
conformal and not Killing, i.e., $f_Y\ne 0$.
 Note that $Y$ is a transversal conformal field if and only if
\begin{align}\label{eq2-6}
g_Q(\nabla_X\overline Y,Z)+g_Q(\nabla_Z\overline Y,X)=2f_Y
g_Q(X,Z) \quad X,Z\in N\mathcal F.
\end{align}
 Let $\{E_a\}$ be a local orthonormal basic frame of
$N\mathcal F$. Then we have the following lemma.
\begin{lemma} $([\ref{JJ}])$   Let $(M,g_M,\mathcal F)$ be a Riemannian manifold with a foliation $\mathcal F$ of
codimension $q$ and a bundle-like metric $g_M$. If $\overline Y
\in \overline V(\mathcal F)$ is a transversal conformal field,
i.e., $\theta(Y)g_Q=2f_Y g_Q$, then we have
\begin{gather}
(\theta(Y) Ric^\nabla)(E_a,E_b) = -(q-2)\nabla_a f_b + \delta_a^b
(\Delta_B f_Y - \kappa^\sharp(f_Y)),\label{eq4-6}
\end{gather}
where $\nabla_a=\nabla_{E_a}$, $f_a=\nabla_a f_Y$ and
$Ric^\nabla(X,Y)=g_Q(\rho^\nabla(X),Y)$ for any $X,Y\in N\mathcal
F$.
\end{lemma}

\section{The generalized Obata theorem}
Recall the following definition similar to that in [\ref{ken2}],
which is a
special case of an isometric equivalence between two pseudogroups
of local isometries acting on smooth manifolds. Let
$(M,g_M,\mathcal F)$ be a Riemannian manifold of a foliation
$\mathcal F$ and a bundle-like metric $g_M$.

\begin{defn}{\rm Let $G$ be a discrete group. A Riemannian foliation
$\left(M,\mathcal F\right)$ is}  transversally isometric {\rm to
$(W,G)$, where $G$ acts by isometries on a Riemannian manifold
$(W,g_W)$,
 if there
exists a homeomorphism $\eta:W\slash G\to M\slash\mathcal{F}$ that
is} locally covered by isometries{\rm . That is, given any $x\in
M$, there exists a local smooth transversal $V$ containing $x$ and
a neighborhood $U$ in $W$ and an isometry $\phi:U\to V$ such that
the following diagram commutes
\[
\begin{array}{ccc}
U & \overset{\phi}{\longrightarrow} & V \\
\ \ \ \ \downarrow^{P\circ i} & \circlearrowleft & \ \ \ \ \downarrow^{\tilde P\circ j}\\
W\slash G & \overset{\eta}{\longrightarrow} & M\slash \mathcal{F}
\end{array}
\]
where $i:U\to W$ and $j:V\to M$ are inclusions and $P:W\to W\slash
G$ and $\tilde P:M\to M\slash \mathcal{F}$ are the projections.}
\end{defn}
Then we have the following generalized Obata theorem for
foliations (recall that $\nabla$ refers to the connection on
$N\mathcal{F}$).
\begin{thm} Let $(M,g_M,\mathcal F)$ be a connected, complete Riemannian manifold with a
foliation $\mathcal F$ of codimension $q\geq 2$ and a bundle-like
metric $g_M$, and let $c$ be a positive real number. Then the
following are equivalent:

$(1)$ There exists a non-constant basic function $f$ such that
$\nabla_X df=-c^2f X^b$ for all vectors $X\in N(\mathcal F)$.

$(2)$ $(M,\mathcal F)$ is transversally isometric to
$\left(S^{q}\left({1}/{c}\right),G\right)$, where the discrete
subgroup $G$ of the orthogonal group $O(q)$ acts by isometries on
the last $q$ coordinates of the $q$-sphere $S^q(1/c)$ of radius
$1/c$ in Euclidean space $\mathbb R^{q+1}$.
\end{thm}
{\bf Proof.} It is clear that the second condition implies the
first, because if $f$ is the first coordinate function in $\mathbb
R^{q+1}$ considered as a function on the sphere $S^q(1/c)$, it
satisfies the first condition. Next, assume that the first
condition is satisfied for the basic function $f$. This implies
that for each $x\in M$,
\begin{align*}
-c^2 f(x)g_{N_x}=\nabla^2 f|_{N_x\mathcal F},
\end{align*}
where $N_x\mathcal F$ is the normal space to the leaf through
$x\in M$, and where $g_{N_x}=g_{N\mathcal F}|_{N_x\mathcal F}$ is
the metric restricted to $N_x\mathcal F$. For any unit speed
geodesic $\gamma:[0,\beta)\to M$ that is normal to the leaves of
the foliation,
\begin{align*}
-c^2(f\circ\gamma)&=-c^2(f\circ\gamma)g_M(\gamma',\gamma')\\
&=g_M(\nabla_{\gamma'}\mathrm{grad}\, f,\gamma')\\
&=g_M(\mathrm{grad}\, f,\gamma')'-g_M(
\mathrm{grad}\, f,\nabla_{\gamma'}\gamma')\\
&=(f\circ\gamma)''.
\end{align*}
Note that since the metric is bundle-like,
every geodesic with initial velocity in $N\mathcal{F}$
is guaranteed to be orthogonal to $\mathcal{F}$ at
all points (\cite{Rein}).
Thus
\begin{align*}
(f\circ\gamma)(t)=A \cos (ct) +B \sin (ct)
\end{align*}
for some constants $A$ and $B$. Let $\gamma(0)=x_0\in M$ be either
a global maximum or global minimum of $f$ on $M$. Then
\begin{align}\label{eq2-1}
f(\gamma(t))=f(x_0)\cos (ct)
\end{align}
for any unit speed geodesic $\gamma$ orthogonal to the leaf ${
L}_{x_0}$ through $x_0$,
and the maximum and minimum values along $\gamma$ must have opposite signs.
Suppose that we choose
the geodesic so that it connects an absolute maximum $x_0$ with an
absolute minimum $x_1$; such a normal geodesic can always be found
(see [\ref{ken1}]).
Note that the nondegeneracy of the normal
Hessian implies that each maximum and minimum of $f\circ\gamma$ occurs at an
isolated closed leaf of $(M,\mathcal F)$; then
the set $f^{-1}(-f(x_0))$ must be a discrete union of closed leaves. The
normal exponential map is surjective ([\ref{ken1}]), and
$f^{-1}([f(x_0),-f(x_0)])=M$ by the reasoning above. So
$f^{-1}(-f(x_0))$ is a single closed leaf, say ${L}_{x_1}$, so
that all normal geodesics through $x_0$ meet ${L}_{x_1}$ at the
exact distance $\pi/ c$. Similarly, $f^{-1}(f(x_0))={L}_{x_0}$.

Given any leaf  ${L}$ of $M$ that is neither $ {L}_{x_0}$ nor
${L}_{x_1}$, there exists a minimal normal geodesic connecting it
to ${L}_{x_0}$ by completeness. In fact, there exists such a
minimal normal geodesic through $x_0$, and its initial velocity
lies in $N_{x_0}\mathcal F$. By equation (\ref{eq2-1}), the
gradient of $f$ is nonzero at each $\gamma(t)$ for $0<t<{\pi/ c}$
and is parallel to $\gamma'(t)$.  Since geodesics are determined
by velocity at a single point, it is impossible that two geodesics
with initial velocities through $x_0$ meet at the same point
unless that point has distance at least $\pi/ c$ from $x_0$. Thus,
the normal exponential map $\exp _{x_0}^\perp : N_{x_0}\mathcal
F\to M$ is injective on the ball $B_{\pi/ c}:=B_{\pi/
c}(x_0)\subset N_{x_0}\mathcal F$. This discussion is independent
of the initial point of ${L}_{x_0}$ chosen, because for a
bundle-like metric the distance from a point $x_0$ on one leaf
closure to another is independent of the choice $x_0\in L_{x_0}$
(see [\ref{ken1}]). We have $\cup_{x\in L_{x_0}} \exp_{x}^\perp
(\overline {B_{\pi/c}(x)})=M$. Let $M_s=\{L_y |\
dist(L_{x_0},L_y)=s\}$ for any non-negative real number $s$, so
that $M_0={L}_{x_0}$ and $M_{\pi/c}={L}_{x_1}$. By the preceding
discussion, for $s\in (0,\pi/c)$, $M_s$ is diffeomorphic to the
unit normal sphere bundle of ${L}_{x_0}\subset M$. Note that the
infinitessimal holonomy group $G$ at $x_0$ acts by orthogonal
transformations on $N_{x_0}\mathcal F$ (\cite{Molino}), and this
action induces an isometric group action on $M_s\cap
\exp_{x_0}^\perp(N_{x_0}\mathcal F)$, with the induced metric from
$g_{N\mathcal{F}}$. Each saturated submanifold $M_s$ for $0\leq s
<{\pi}/{c}$ has a leaf space that is isometric to the quotient of
$M_s\cap \exp_{x_0}^\perp(N_{x_0}\mathcal F)$ by $G$. Then
$(M\smallsetminus L_{x_1})/\mathcal F$ is diffeomorphic to
$B_{\pi/c}/G$. The map
\begin{align*}
\eta : B_{\pi/c}/G \to (M\smallsetminus L_{x_1})/\mathcal F
\end{align*}
is defined by $\eta(O_\xi)=L_{exp_{x_0}^\perp(\xi)}$, where
$\xi\in B_{\pi/c}\subset N_{x_0}\mathcal F$, $O_\xi$ is the
$G$-orbit of $\xi$ in $B_{\pi/c}$ and $L_{exp_{x_0}^\perp(\xi)}$
is the leaf containing  $\exp_{x_0}^\perp(\xi)$. Letting $B_{\pi/c}^{+}$
denote the one-point compactification of $B_{\pi/c}$,
$\eta$ can be extended to a homeomorphism
\begin{align*}
\overline \eta : B_{\pi/c}^+/G \to M/\mathcal F.
\end{align*}
Thus $M/\mathcal F$ is homeomorphic to $S/G$, where
$S=B_{\pi/c}^+$ is a sphere. Next we will show that the pullback
of the transverse metric of $\left(M,\mathcal{F}\right)$ endows
$S$ with the standard metric of $S^q\left({1}/{c}\right)$.

Let $v$ and $w$ be any two nonzero orthonormal vectors in
$N_{x_0}\mathcal F$, and let $W_s$ denote the
$N\mathcal{F}$-parallel translate of $w=W_0$ along the geodesic
$\gamma(s)$ with initial velocity $v$; thus $W_s\in
N_{\gamma(s)}\mathcal{F}$ is a well-defined vector at each
$\gamma(s)$ for $0\leq s<{\pi}/{c}$. We see that $W_s$ is tangent
to $M_s$ for $s\in (0,\pi/c)$. Let $\left( y_j\right)$ be geodesic
normal coordinates for the normal ball $\exp^\perp_{x_0}\left(
B_{\pi/c}\left(x_0\right)\right)$. Suppose that these coordinates
are chosen at $x_0$ such that $y_1(\gamma(s))=s$ and each of
$\frac{\partial}{\partial y_j}$ for $j>1$ is orthogonal to
$v=\gamma'(0)$ at $x_0=0$. We extend $s$ to be the function
$s(y)=\sqrt{\sum y_j^2}$ and write $y_j=s\theta_j$, so that each
$\theta_j$ is independent of $s$. Thus, $\gamma'(s)(\theta_j)=0$
and $W_s(s)=0$. Further, we let $\frac{\partial}{\partial s}$
denote the radial vector field, which agrees with $\gamma'(s)$
along $\gamma$. In the calculations that follow, we extend $y_j$,
$\theta_j$, $\frac{\partial}{\partial s}$ to be well-defined and
basic in a small neighborhood of the transversal
$\exp^\perp_{x_0}\left( B_{\pi/c}\right)$. From the calculation of
$f$ above, we see that $\mathrm{grad}\, f=-c\sin(cs)
f(x_0)\frac{\partial}{\partial s}$.

Since $\nabla$ is torsion-free and $\nabla_{\gamma'(s)}W_s=0$ by
construction,
\begin{align*}
\pi\left[\frac{\partial}{\partial s},W_s\right]&=-\nabla_{W_s}\tfrac{\partial}{\partial s}=\frac 1{c\sin(cs) f(x_0)}\nabla_{W_s}\mathrm{grad}\, f\\
&=-\frac{c^2}{c\sin(cs)f(x_0)}f(\gamma(s))W_s\\
&=-\frac{c\cos(cs)}{\sin(cs)}W_s,
\end{align*}
by assumption, since $\nabla$ is a metric connection and thus
commutes with raising indices. Since $\theta_j$ is a locally
defined basic function, for $0<s<{\pi}/{c}$,
\begin{align*}
\frac{d}{ds}W_s\left(\theta_j\right)=\frac{\partial}{\partial
s}W_s\left(\theta_j\right) =\left[\frac{\partial}{\partial
s},W_s\right](\theta_j) =\pi\left[\frac{\partial}{\partial
s},W_s\right](\theta_j)
=-\frac{c\cos(cs)}{\sin(cs)}W_s\left(\theta_j\right).
\end{align*}
Solving the differential equation above, we have
\begin{align}
W_s\left(\theta_j\right)=\frac{1}{\sin cs}W_{\pi/2c}\left(\theta_j\right),
0<s<\frac\pi c.
\end{align}
Since $W_s(s)=0$ we have
\[
W_s\left(y_j\right)=sW_s\left(\theta_j\right)
\]
for $0<s<\pi/ c$. Then, for all $j$,
\begin{align*}
W_0\left(y_j\right)
&=\lim_{s\to 0}W_s\left(y_j\right) \\
&= \frac 1c W_{\frac{\pi}{2c}} \left(\theta_j\right)\\
&= \frac {\sin (cs)}{c} W_s\left(\theta_j\right)\\
&=\frac {\sin (cs)}{cs} W_s\left(y_j\right).
\end{align*}
Note that since the vectors $\frac{\partial}{\partial \theta_j}$
for $j>1$ form a basis of the tangent space for
$M_s\cap\exp^\perp_{x_0}\left(B_{\pi/c}\right)$ at $\gamma(s)$
with $s>0$, the equation above uniquely defines the vector $W_s$
in terms of $W_0$. Since the metric on the sphere
$S^q\left({1}/{c}\right)$ satisfies the same hypothesis, a
corresponding fact is true for geodesic normal coordinates on
$S^q\left({1}/{c}\right)$.

We now show that the equation above implies that the pullback of
the metric $g_{N\mathcal{F}}$ to $B_{\pi/c}$ is the same as the
standard metric $g_S$ corresponding to geodesic normal coordinates
on $S^q\left({1}/{c}\right)$. As above, let $W_s$ denote the
parallel displacement of $W_0$ along $\gamma(s)$, and let
$\overline{W_s}$ denote the parallel displacement of ${W_0}$ along
the geodesic in $\left(B_{\pi/c},g_S\right)$ with unit tangent
vector $v$. Then
\[
W_s\left({y_j}\right)
=\frac{cs}{\sin(cs)}W_0\left(y_j\right)
=\overline{W_s}\left({y_j}\right).
\]
Since the actions of the vectors $\overline{W_s}$ on the
coordinate functions ${y_j}$ determine their values along the
geodesic with initial velocity $v$, we conclude that
$\overline{W_s}=W_s$. Thus, the metrics $g_{N\mathcal{F}}$ and
$g_S$ on $B_{\pi/c}$ yield identical parallel displacements of
vectors orthogonal to $v$ along the line containing $v$, and
$\left.g_{N\mathcal{F}}\right|_{x_0}=\left.g_{S}\right|_{x_0}$.
Since it follows from previous calculations that the initial
vector $v\in N_{x_0}\mathcal{F}$ is arbitrary, we conclude that
$g_{N\mathcal{F}}=g_S$. We may reverse the roles of $x_0$ and
$x_1$ and obtain a similar result.

Now, given any point $x\in M$, there is a minimal geodesic
connecting this point to a point $x_0'$ on the leaf containing
$x_0$. If $x\notin L_{x_1}$, the above analysis shows that the map
$\exp^\perp_{x_0'}$ restricted to
$\left(B_{\pi/c}\left(x_0'\right),g_S\right)$ is an isometry onto
its image, and that image contains $x$. Further, the map
$\exp^\perp_{x_0'}$ locally covers the map
$\overline\eta:\left(B_{\pi/c}^+\slash G,g_S\right) \to
\left(M\slash\mathcal F,g_{N\mathcal{F}}\right)$. If $x\notin
L_{x_0}$, a similar fact is true for $\exp^\perp_{x_1'}$. Thus the
map $\overline\eta$ is locally covered by isometries, and we
conclude that $(M,\mathcal F)$ is transversally isometric to
$\left(S^{q}\left({1}/{c}\right),G\right)$. $\Box$

\section{Applications}
In this section, we give some applications of the generalized
Obata theorem. Let $M$ be a Riemannian manifold admitting a
transversal non-isometric conformal field. For more details about
transversal conformal fields, see [\ref{JJ},\ref{Pak}].

\begin{thm}\label {thm3-5}  Let $(M,g_M,\mathcal F)$ be a compact Riemannian manifold with a foliation $\mathcal F$ of
codimension $q$ and a bundle-like metric $g_M$. Assume that the
transversal scalar curvature $\sigma^\nabla$ is a positive
constant. If $M$ admits a transversal non-isometric conformal
field $\overline Y$ such that $\overline Y=\nabla h$ for some
basic function $h$, then $(M,\mathcal F)$ is transversally
isometric to $(S^q(1/c),G)$, where $c^2={\sigma^\nabla\over
q(q-1)}$ and $G$ is a discrete subgroup of $O(q)$ acting on the
$q$-sphere .
\end{thm}
{\bf Proof.} Let $\theta(Y) g_Q=2f_Y g_Q (f_Y\ne 0)$  and
$\overline Y=\nabla h$ for some basic function $h$. From
(\ref{eq2-6}), we have that, for any $X,Z\in\Gamma Q$,
\begin{align*}
2f_Y g_Q(X,Z)&=g_Q(\nabla_X\nabla h,Z)+ g_Q(\nabla_Z\nabla h,X)\\
&=\nabla\nabla h(X,Z)+\nabla\nabla h(Z,X)\\
&=2\nabla\nabla h(X,Z).
\end{align*}
Hence we have
\begin{align}\label{eq4-8}
\nabla\nabla h= f_Y g_Q.
\end{align}
Note that the function $f_Y$ satisfies [\ref{JJ}]
\begin{align}\label{eq4-7}
\Delta_B f_Y = {\sigma^\nabla \over q-1}f_Y +\kappa_B^\sharp(f_Y).
\end{align}
Since $\Delta_B h=-\sum_a\nabla_{E_a}\nabla_{E_a} h
+\kappa_B^\sharp(h)$[\ref{JJ}], we have  from (\ref{eq4-8}) and
(\ref{eq4-7})
\begin{align*}
\Delta_B\Big(f_Y+{\sigma^\nabla\over
q(q-1)}h\Big)=\kappa_B^\sharp\Big(f_Y +{\sigma^\nabla\over
q(q-1)}h\Big).
\end{align*}
From Lemma 2.1, we have that
\begin{align}
f_Y + {\sigma^\nabla\over q(q-1)}h={\rm constant},
\end{align}
which yields
\begin{align*}
\nabla\nabla f_Y +{\sigma^\nabla\over q(q-1)}\nabla\nabla h=0.
\end{align*}
From (\ref{eq4-8}), we have
\begin{align}
\nabla\nabla f_Y = -{\sigma^\nabla\over q(q-1)}f_Y g_Q.
\end{align}
By the generalized Obata theorem(Theorem 3.2), $(M,\mathcal F)$ is
transversally isometric to $(S^q(1/c),G)$, where
$c^2={\sigma^\nabla\over q(q-1)}$. $\Box$

\begin{thm}  Let $(M,g_M,\mathcal F)$ be a compact Riemannian manifold with a foliation $\mathcal F$ of
codimension $q$ and a bundle-like metric $g_M$. Assume that the
transversal scalar curvature $\sigma^\nabla$ is a positive
constant. If $M$ admits a transversal non-isometric conformal
field $\overline Y$ such that $\theta(Y) Ric^\nabla = \mu g_Q$ for
some basic function $\mu$, then $(M,\mathcal F)$ is transversally
isometric to $(S^q(1/c),G)$, where $c^2={\sigma^\nabla\over
q(q-1)}$ and $G$ is a discrete subgroup of $O(q)$ acting on the
$q$-sphere.
\end{thm}
{\bf Proof.} Let $\theta(Y)g_Q=2f_Y g_Q (f_Y\ne 0)$. From Lemma
2.2, we have
\begin{align}\label{eq4-11}
\mu g_Q = -(q-2)\nabla\nabla f_Y + (\Delta_B f_Y
-\kappa_B^\sharp(f_Y))g_Q.
\end{align}
Hence we have
\begin{align}\label{eq4-12}
\mu = {2(q-1)\over q} (\Delta_B f_Y-\kappa_B^\sharp(f_Y)).
\end{align}
From (\ref{eq4-11}) and (\ref{eq4-12}), we have
\begin{align}\label{eq4-13}
\nabla\nabla f_Y = -\frac1q (\Delta_B f_Y
-\kappa_B^\sharp(f_Y))g_Q.
\end{align}
From (\ref{eq4-7}), we have
\begin{align}
\nabla\nabla f_Y + {\sigma^\nabla\over q(q-1)}f_Y g_Q=0.
\end{align}
By the generalized Obata theorem, the proof is completed.  $\Box$

On the other hand, we recall the following theorem.
\begin{thm}([\ref{JJ}]) Let $(M,g_M,\mathcal F)$ be a compact Riemannian
manifold with a foliation $\mathcal F$ of codimension $q$ and a
bundle-like metric $g_M$ with $\delta_B\kappa_B=0$. Assume that
the transversal scalar curvature $\sigma^\nabla$ is constant and
$\rho^\nabla(X)\geq {\sigma^\nabla\over q}X$ for any $X\in \Gamma
Q$. If $M$ admits a transversal non-isometric conformal field,
then $(M,\mathcal F)$ is transversally isometric to $(S^q,G)$,
where $G$ is a discrete subgroup of $O(q)$ acting on the
$q$-sphere.
\end{thm}
{\bf Remark.} On a compact Riemannian manifold admitting a
transversal non-isometric conformal field, if the scalar curvature
$\sigma^\nabla$ is constant, then the condition
$\delta_B\kappa_B=0$ implies that $\sigma^\nabla$ is non-negative.
Moreover, if $\sigma^\nabla$ is positive constant and the
transversal Ricci curvature $\rho^\nabla(X)\geq
{\sigma^\nabla\over q}X$ for any $X\in \Gamma Q$, then $\kappa=0$
by the tautness theorem ([\ref{Min}]).

Hence Theorem 4.3 is equivalent to the following.
\begin{thm}Let $(M,g_M,\mathcal F)$ be a compact Riemannian
manifold with a foliation $\mathcal F$ of codimension $q$ and a
bundle-like metric $g_M$. Assume that the transversal scalar
curvature $\sigma^\nabla$ is a positive constant and
$\rho^\nabla(X)\geq {\sigma^\nabla\over q}X$ for any $X\in \Gamma
Q$. If $M$ admits a transversal non-isometric conformal field,
then $(M,\mathcal F)$ is transversally isometric to $(S^q,G)$,
where $G$ is a discrete subgroup of $O(q)$ acting on the
$q$-sphere.
\end{thm}
We define an operator $A_Y:N\mathcal F \to N\mathcal F$ for any
vector field $Y\in V(\mathcal F)$ by
\begin{equation}
 A_Y s=\theta(Y)s -\nabla_Y s.
\end{equation}
  Then it is proved [\ref{Kamber2}] that, for any vector field
$Y\in V(\mathcal F)$,
\begin{equation}
A_Y s = -\nabla_{Y_s}\bar Y,
\end{equation}
where $Y_s =\sigma(s)\in \Gamma TM$.   So $A_Y$ depends only on
$\bar Y=\pi(Y)$ and is a linear operator. Moreover, $A_Y$ extends
in an obvious way to tensors of any type on $N\mathcal F$ (see
[\ref{Kamber2}] for details).
\begin{thm} Let $(M,g_M,\mathcal F)$ be a compact Riemannian
manifold with a foliation $\mathcal F$ of codimension $q$ and a
bundle-like metric $g_M$. Assume that the transversal scalar
curvature $\sigma^\nabla$ is a positive constant. If $M$ admits a
transversal non-isometric conformal field $\overline Y$, i.e.,
$\theta(Y)g_Q=2f_Y g_Q$ such that (i) $\rho^\nabla(\nabla
f_Y)={\sigma^\nabla\over q}\nabla f_Y$, (ii) $\kappa_B(f_Y)=0$ and
(iii) $g_Q(A_{\kappa_B^\sharp}\nabla f_Y, \nabla f_Y)\leq 0$, then
$(M,\mathcal F)$ is transversally isometric to $(S^q,G)$, where
$G$ is a discrete subgroup of $O(q)$ acting on the $q$-sphere.
\end{thm}
{\bf Proof.} First, we recall that, for any basic function $g$
with $\kappa_B^\sharp(g)=0$, we have [\ref{JJ}]
\begin{align*}
\int_M\{&{q-1\over q}g_Q(\Delta_B d_B g, d_B
g)-g_Q(\rho^\nabla(\nabla f_Y),\nabla
f_Y)+g_Q(A_{\kappa_B^\sharp}\nabla f_Y, \nabla f_Y)\\
&-|\nabla\nabla g +{1\over q}\Delta_B g |^2\}=0.
\end{align*}
From (\ref{eq4-7}) and assumption (ii), we have
\begin{align*}
\Delta_B f_Y = {\sigma^\nabla \over q-1}f_Y.
\end{align*}
Hence, from the assumptions (i) and (iii), we have
\begin{align*}
\nabla\nabla f_Y =-{\sigma^\nabla\over q(q-1)}f_Y g_Q.
\end{align*}
By the generalized Obata theorem, the proof is completed. $\Box$

If the foliation $\mathcal F$ is minimal, the conditions $(ii)$
and $(iii)$ in Theorem 4.5 are satisfied. Hence we have the
following corollaries.
\begin{coro} Let $(M,g_M,\mathcal F)$ be a compact Riemannian
manifold with a minimal foliation $\mathcal F$ of codimension $q$
and a bundle-like metric $g_M$. Assume that the transversal scalar
curvature $\sigma^\nabla$ is a positive constant. If $M$ admits a
transversal conformal field $\overline Y$, i.e.,
$\theta(Y)g_Q=2f_Y g_Q$ such that $\rho^\nabla(\nabla
f_Y)={\sigma^\nabla\over q}\nabla f_Y$, then $(M,\mathcal F)$ is
transversally isometric to $(S^q,G)$, where $G$  a discrete
subgroup of $O(q)$ acting on the $q$-sphere.
\end{coro}

\begin{coro} Let $(M,g_M,\mathcal F)$ be a compact Riemannian
manifold with a minimal, transversally Einstein foliation
$\mathcal F$ of codimension $q$ and a bundle-like metric $g_M$.
Assume that the transversal scalar curvature $\sigma^\nabla$ is a
positive constant. If $M$ admits a transversal non-isometric
conformal field $\overline Y$, then $\mathcal F$ is transversally
isometric to $(S^q,G)$, where $G$ is a discrete subgroup of $O(q)$
acting on the $q$-sphere.
\end{coro}

\noindent Department of Mathematics, Jeju National University,
Jeju 690-756, Korea

\noindent {\it E-mail address} : sdjung@jejunu.ac.kr

\noindent Department of Mathematics, Jeju National University,
Jeju 690-756, Korea

\noindent{\it E-mail address} :niver486@jejunu.ac.kr

\noindent Department of Mathematics, Texas Christian University,
Fort Worth, TX 76129, USA

\noindent{\it E-mail address} : k.richardson@tcu.edu

\end{document}